\documentclass[12pt]{amsart}
\usepackage{amsfonts,amsmath,amstext,amsthm,amssymb}
\usepackage{hyperref}
\numberwithin{equation}{section}
\newtheorem{theorem}{Theorem}
\newtheorem*{theo}{Theorem}
\newtheorem{lemma}{Lemma}

\def\sign{{\rm sign\,}}
\def\diam{{\rm diam\,}}

\def\supp{{\rm supp\,}}
\def\rot{{\rm rot\,}}
\def\dil{{\rm dil \,}}

\def\mes{{\rm mes}}
\def\tg{{\rm tg }}
\def\MR{\ensuremath{\mathcal R}}
\def\ZR{\ensuremath{\mathbb R}}
\def\ZZ{\ensuremath{\mathbb Z}}

\def\ZN{\ensuremath{\mathbb N}}

\def\ZI{\ensuremath{\mathbb I}}

\newcommand {\e }[1]{(\ref{#1})}
\newcommand {\lem }[1]{Lemma \ref{#1}}

\begin{document}
 \baselineskip=17pt

\title[A complete characterization of
$R$-sets]{A complete characterization of \\
$R$-sets in the theory of differentiation of integrals}

\author
{G. A. Karagulyan}

\address{Yerevan State University, Depart. of Computer Science}

\address{Institute of Mathematics \\
Armenian National Academy of Sciences\\ Marshal Baghramian ave.
24b,\\ Yerevan, 375019, ARMENIA\\ }

\email{karagul@instmath.sci.am}
\date{}
\maketitle
\begin{abstract}
Let $\MR_s$ be the family of open rectangles in the plane $R^2$
having slope $s$ with the abscissa. We say a set of slopes $S$ is
$R$-set if there exists a function $f\in L(R^2)$, such that the
basis $\MR_s$ differentiates integral of $f$ if $s\not\in S $, and
$$
\overline D_sf(x)=\limsup_{\diam(R)\to 0, x\in R\in\MR_s}
\frac{1}{|R|}\int_R f=\infty
$$
almost everywhere if $s\in S$. If
the condition $\overline D_s f(x)=\infty $ holds on a set of
positive measure (instead of a.e.) we shall say it is $WR$-set. It
is proved, that $S $ is a $R$-set($WR$-set) if and only if it is
$G_\delta $($G_{\delta\sigma}$).
\end{abstract}

\begin{section}{Introduction}\everypar{} \parskip=11pt
For any number $s\in [0,\frac{\pi}{2})$ we define $\MR_s$ to be
the family of all open rectangles $R$ in $\ZR^2$ having slope $s$,
i.e. $R$ has a side forming angle $s$ with the abscissa. We say
that the basis $\MR_s$ differentiates the integral of the function
$f\in L^1(\ZR^2 )$, if
\begin{equation}\label{fdif}
 \lim\limits_{d(R)\to 0, x\in R\in\MR_s} \frac{1}{|R|}\int_R
f=f(x)
\end{equation}
almost everywhere in $\ZR^2 $, where $d(R)$ is the diameter of
$R$. According to the well-known theorem of
Jessen-Marcinkiewicz-Zygmund \cite{JMZ} the basis $\MR_s$
differentiates $\int f$ for any function $f\in L\log L(\ZR^2)$. On
the other hand S.~Saks \cite{Saks} constructed an example of
function $f\in L^1(\ZR^2 )$ such that
\begin{equation*}
\overline D_sf(x)=\limsup_{d(R)\to 0, x\in R\in\MR_s}
\frac{1}{|R|}\int_R f=\infty ,\hbox { everywhere }.
\end{equation*}
In view of this A.~Zygmund in \cite{Guz} posed the following
problem: for a given $f\in L^1(\ZR^2 )$ is it possible to find a
direction $s$ such that $\MR_s$ differentiates $\int f$?
J.~Marstrand in \cite{Mar} gave a negative answer to this
question, proving
\begin{theo}[J.~Marstrand]
There exists a function $f\in L^1(\ZR^2 )$ such that $\overline
D_sf(x)=\infty $ almost everywhere for any $s$.
\end{theo}
\noindent
 Different generalizations of this result are obtained by J.~El Helou
\cite{Hel}, A.~M.~Stokolos \cite{Sto}, B.~L{\'o}pez Melero
\cite{Mel} and G.~G.~Oniani \cite{Oni2}. A.~M.~Stokolos in
\cite{Sto} extended Marstrand's theorem to higher dimensional
case. In the papers \cite{Mel} and \cite{Oni2} it is considered
the same problem for general translation invariant differentiation
basises.

\noindent We say that the set $S\subset [0,\frac{\pi}{2})$ is
$R$-set if there exists a function $f\in L^1(\ZR^2 )$ such that
the basis $\MR_s$ differentiates  $\int f$ whenever $s\in [0,\pi
/2)\setminus S$, and $\overline D_sf(x)=\infty $ almost everywhere
as $s\in S$. If the condition $\overline D_sf(x)=\infty $ holds on
a set of positive measure (instead of a.e.) we shall say it is
$WR$-set (weak $R$-set).
 In this language, Marstrand's theorem asserts,
that $[0,\pi/2)$ is $R$-set. A.~M.~Stokolos in \cite{Sto2} proved,
the existence of everywhere dense $WR$-set, which is not whole
$[0,\pi/2)$. G.~Lepsveridze in \cite{Lep1},\cite{Lep2}  proved
that any finite set is $R$-set and any countable set is in some
$WR$-set of measure zero. G.~G.~Oniani in \cite{Oni2} generalizing
this result proved that any countable set is in some $R$-set of
measure zero.

\noindent
 The definition of $R$-sets first appeared in the paper \cite{Oni1} by
 G.~G.~Oniani, where the author posed the problem about characterization of
 all $R$-sets. In particular, it was a question if there exists a
$R$-set of positive measure and moreover whether any interval is
$R$-set or not? In the same paper Oniani shows, that any $R$-set
is $G_\delta $ in $[0,\pi/2)$, i.e.
\begin{equation*}
G=\big(\cap_{k=1}^\infty G_k\big)\cap [0,\pi/2)
\end{equation*}
where $G_n$ are open sets, and conversely if $G_\delta $-set is
countable, then it is $R$-set. These results characterize the
countable $R$-sets. We note that any countable $G_\delta $-set is
nowhere dense. So in \cite{Oni1} Oniani constructed also a $R$-set
of second category. These problems are stated also in the
monograph G.~G.~Oniani \cite{Oni2} and in the papers \cite{Oni3}
and \cite{Oni4} it is investigated the higher dimensional case of
the problem.

 \noindent The following theorems give a complete characterization of
 general $R$ and $WR$-sets.
 \begin{theorem}\label{Th} For the set $S\subset [0,\pi/2 )$ to be $R$-set
 it is necessary and sufficient to be $G_\delta $.
 \end{theorem}
 \begin{theorem}\label{Th} For the set $S\subset [0,\pi/2 )$ to be $WR$-set
 it is necessary and sufficient to be $G_{\delta\sigma} $.
 \end{theorem}

\noindent The necessity of Theorem 1 is proved by Oniani in
\cite{Oni1}. We present here a short statement of the proof of
that. If $S$ is a $R$-set, then there exists a function $f\in L^1$
such that \e{fdif} holds as $s\in [0,\pi/2)\setminus S$ and
$\overline D_sf(x)=\infty $ a.e. as $s\in S$. For any $n\in \ZN $
denote
\begin{equation*}
U_n=\{s\in [0,\pi/2):\, |\{x\in B(n),\,
M_s^{[0,1/n)}f(x)>n\}|>|B(n)|-2^{-n}\},
\end{equation*}
where $B(n)=\{x\in \ZR^2:\|x\|\le n\}$ and the maximal function
$M_sf$ are defined in Section 2. It is easy to check, that
$U_n=G_n\cap [0,\pi/2)$, where $G_n$ are open sets and
\begin{equation*}
\{s\in [0,\pi/2):\,\overline D_sf(x)=\infty\, \text { a.e. }
\}=\bigcap_nU_n=\bigg(\bigcap_nG_n\bigg)\bigcap [0,\pi/2) ,
\end{equation*}
i.e. it is $G_\delta$-set in $[0,\pi/2)$, which proves the one
part of Theorem 1.

To prove the necessity of Theorem 2 it is enough to prove that for
any function $f\in L^1(\ZR^2)$ the set
\begin{equation*}
G_f=\{s\in [0,\pi/2):\,|\{x\in \ZR^2:\, \overline D_sf(x)=\infty
\}|>0 \}
\end{equation*}
is $G_{\delta\sigma}$. Denote
\begin{equation*}
U_{nm}=\{s\in [0,\pi/2):\,|\{x\in B(n):\, M_s^{[0,1/m)}f(x)>m \}|>\frac{1}{n} \},\quad
n,m=1,2,\cdots ,
\end{equation*}
where $B(n)$ and $M_sf$ are defined in Section 2. It is clear
$U_{nm}$ are open sets in $[0,\pi/2 )$ and
\begin{equation*}
G_f=\bigcup_n\bigcap_m U_{nm}.
\end{equation*}
To show the last equality it suffices to check the following
relations:
\begin{multline*}
s\in G_f\Leftrightarrow |\{x\in \ZR^2:\, \overline D_sf(x)=\infty
\}|>\alpha >0\\ \Leftrightarrow \exists\, n\hbox { such that }
|\{x\in B_n:\, \overline D_sf(x)=\infty \}|>\frac{1}{n}\\
\Leftrightarrow \exists\, n\hbox { such that }s\in
\bigcap_mU_{n,m} \Leftrightarrow s\in \bigcup_n\bigcap_m U_{n,m}.
\end{multline*}
Hence the set $G_f$ is $G_{\delta\sigma }$.

 We shall prove the
sufficiencies of the theorems invoking the probabilistically
independence of sets similar to original approach of J.~Marstrand
in \cite{Mar}. This idea is involved in \lem{lm5}. Of coarse, we
use also Bohr's construction displayed in Saks's classical
counterexample. It is important that the function constructed  in
the proof is not nonnegative, which we don't have in all the
results stated above. This argument gives more freedom in the
construction to ensure differentiability of the integral along
some directions. So the method demonstrated in the proof differs
from the others, because we essentially use an interference of
positive and negative values of a function in integrals, which is
displayed in \lem{lmR} and \lem{lm2}.

\end{section}
\begin{section}{Notations and Lemmas}
\noindent The basis $\MR_s$ can be defined for any  $s\in [0,2\pi
]$. We note that $\MR_s=\MR_t$ if $s=t\mod\pi/2 $. In fact
$\cup_{s\in[0,\pi/2)}\MR_s$ is the family of all rectangles in the
plane.

\noindent If $n\in \ZN$ is an integer and $c=(c_1,c_2)$, then for
any set $A\subset \ZR^2$ we denote
\begin{gather*}
\dil_nA=\{x=(x_1,x_2)\in \ZR^2:\, nx=(nx_1,nx_2)\in A\},
\\
c+A=\{x=(x_1,x_2)\in \ZR^2:\, x=c+a,\, a\in A\}.
\end{gather*}
\noindent We let $Q_0=[-1/2,1/2)\times [-1/2,1/2)$ and for any
$n\in \ZN$, $k=(k_1,k_2)\in\ZZ^2$ denote
$Q_k^n=\dil_n\big(k+Q_0\big)$. For a fixed $n$ the family
$\{Q_k^n: k\in\ZZ^2\}$ is a partition of the plane to squares with
side lengths $1/n$. In some places for $Q_k^1$ we shall use simply
$Q_k$.

\noindent We denote by $\rot_sA$ the rotation of the set $A\subset
\ZR^2$ round the point $(0,0)$ by angle $s$. Denote $B(\varepsilon
)=\{x\in\ZR^2:\, \|x\|=\sqrt{x_1^2+x_2^2}\le\varepsilon \}$ and
$\Gamma_s(\varepsilon )=\rot_s\{x=(x_1,x_2):\,|x_2|<\varepsilon
\}$.

\noindent The notation $s^\perp $ stands for the direction
$s+\pi/2$. For any direction $s$ define $\mes_sA$ to be the linear
Lebesgue measure of the projection of $A$ on the line parallel to
$s^\perp $.

\noindent For any measurable set $A\subset \ZR^2$ we denote
\begin{gather*}
\mes^*A=\sup_{k\in\ZZ^2} |A\cap Q_k|,
\\
\mes_*A=\inf_{k\in\ZZ^2} |A\cap Q_k|.
\end{gather*}

\noindent For numbers $0<\delta<\mu\le\infty $ we define
$\MR_s^{[\delta,\mu )}$ to be the family of rectangles
$R=R_1\times R_2\in\MR_s$ with $\delta\le |R_1|,|R_2|<\mu $ and we
let $\MR_s^\delta$ to be the rectangles from $\MR_s$ with
$|R_1|=|R_2|=\delta $. Denote
\begin{equation*}
 M_s^{[\delta, \mu )}f(x)=\sup_{R\in\MR_s^{[\delta,\mu
)}}\frac{1}{|R|}\bigg|\int_Rf(x)dx\bigg|.
\end{equation*}
If $\delta =0$ and $\mu=\infty $ we shall use notation $M_s f(x)$.
We say that the set $A\subset \ZR^2$ is $\delta $-set if it is a
union of mutually disjoint rectangles from the family
$\MR_s^{[\delta,\infty )}$. The following lemma contains the main
idea of the proof of Marstrand's theorem.
\begin{lemma}\label{lm5}
Suppose $0<\delta_t<1,\, t=1,2,\cdots ,T$ are arbitrary numbers
and $A_t\subset \ZR^2$ are $\delta_t $-sets with
$\mes_*(A_t)>12\delta_t $, $t=1,2,\cdots ,T$. Then for any
sequence of integers $\{n_t\}$, $n_1=1$,
$n_{t+1}>\frac{4}{\delta_t}n_t$, we have
\begin{equation}\label{dilA}
\mes_*\bigg(\bigcup_{t=1}^T \dil_{n_t}
(A_t)\bigg)>1-\bigg(1-\frac{\mes_*(A_t) }{32}\bigg)^T.
\end{equation}
\end{lemma}
\begin{proof}
First we prove that if $B$ is $\delta $-set with $\mes_*B>12\delta
$, $m,n\in \ZN$ and $n>\frac{4}{\delta}m$, then there exists a set
$\tilde B$ such that

\begin{itemize}
    \item [1)]$\tilde B\subset \dil_m B$,
    \item [2)]for any  $k\in \ZZ^2$ the set $\tilde B\cap Q_k^m$ is a union of
    squares $Q_j^n$,
    \item [3)]the values $|\tilde B\cap Q_k^m|$ are equal for
    different $k\in \ZZ^2$,
    \item [4)]$\mes_*(\tilde B)>\frac{1 }{32}\mes_* B$.
\end{itemize}

 We note that any rectangle
$R\in\MR_s^{[\delta,\infty)}$ is a union of rectangles from
$\MR_s^\delta $. So we have $\dil_mB=\cup_i R_i$ where
$R_i\in\MR_s^{\delta /m}$. Denote
\begin{equation*}
B'=\bigcup_{R_i\subset Q_k^m \hbox { for some } k\in
\ZZ^2}R_i\subset \dil_m B.
\end{equation*}
We have $\diam (R_i)=\frac{\delta \sqrt 2}{m}$. So if
$R_i\not\subset Q_k^m$ then $R_i\cap\tilde Q_k^m=\varnothing $ as
$k\in \ZZ^2$, where $\tilde Q_k^m$ is the square concentric
$Q_k^m$ with side lengths  $\frac{1}{m}(1- 2\delta \sqrt 2)$.
Hence we get
\begin{multline}\label{23}
| B'\cap Q_k^m|>|\dil_mB\cap Q_k^m|-|Q_k^m\setminus\tilde Q_k^m|\\
=|\dil_mB \cap Q_k^m|-\frac{1}{m^2}(4\delta\sqrt 2-8\delta^2)
>|\dil_mB \cap Q_k^m|-\frac{6\delta}{m^2}\\
=\frac{1}{m^2}|B \cap Q_k^1|-\frac{6\delta}{m^2} \ge
\frac{1}{m^2}(\mes_*B-6\delta )>\frac{\mes_*B}{2m^2}.
\end{multline}
Using Besicovitch theorem on covering by squares (see \cite{Guz},
p. 10), we may choose a subfamily $\{R'_i\}$ from $\{R_i\}$ such
that $R'_i$ are pairwise disjoin and
\begin{equation}\label{241}
\left|\bigcup_{R'_i\subset Q_k^m }R'_i\right|\ge
\frac{1}{4}\left|\bigcup_{R_i\subset Q_k^m }R_i\right| \hbox{ for
any } k\in \ZZ^2.
\end{equation}
Therefore, denoting
\begin{equation*}
B''=\bigcup R'_i\subset B'\subset \dil_mB,
\end{equation*}
by \e{23} and \e{241} we have
\begin{equation}\label{251}
|B''\cap Q_k^m|>\frac{\mes_*B}{8m^2},\quad k\in \ZZ^2.
\end{equation}
Using a simple geometry, one can easily check that if $R\in
\MR_s^{\delta/m} $ and $n>\frac{4}{\delta }m$, then
\begin{equation*}
\left|\bigcup_{Q_j^n\subset R}Q_j^n\right|>\frac{1}{4}|R|.
\end{equation*}
So, by virtue of \e{251}, for $n>\frac{4}{\delta }m$ we have
\begin{equation*}
\left|\bigcup_{Q_j^n\subset B''\cap
Q_k^m}Q_j^n\right|>\frac{1}{4}|B''\cap
Q_k^m|>\frac{\mes_*B}{32m^2}.
\end{equation*}
Taking away some of the squares $Q_j^n$ from the left union we can
get a set $\tilde B\subset B''$, which is again a union of the
squares $Q_j^n$ and in addition all the sets $\tilde B\cap Q_k^m$
consist of a same number of squares $Q_j^n$ and $|\tilde B\cap
Q_k^m|\ge \frac{\mes_*B}{32m^2},\, k\in \ZZ^2$. Certainly, $\tilde
B$ satisfies the conditions (1)-(4)

\noindent Taking $n=n_{t+1}$, $m=n_t$, $B=A_{n_t}$, $t=1,2,\cdots
,T$ we get sets $\tilde A_t$, $t=1,2,\cdots ,T$, such that
\begin{itemize}
    \item [1)]$\tilde A_t\subset \dil_{n_t}A_t$,
    \item  [2)]$\tilde A_t\cap Q_k^{n_t}$ is a union of squares
    $Q_j^{n_{t+1}}$ for any $k\in \ZZ^2$,
    \item [3)]the values $|\tilde A_t\cap Q_k^{n_t}|$ are equal for
    different $k\in \ZZ^2$,
    \item [4)]$\mes_*(\tilde A_t)>\frac{\mes_*(A_t)}{32}$.
\end{itemize}
From the conditions $2),3)$ it follows that for the fixed $k\in
\ZZ^2$ the sets $\tilde A_t\cap Q_k,\, t=1,2,\cdots ,T$ are
probabilistically independent. Then by 1) and 4)
\begin{multline*}
\mes_*\bigg(\bigcup_{t=1}^T \dil_{n_t} A_t\bigg)\ge
\mes_*\bigg(\bigcup_{t=1}^T\tilde A_t\bigg)=\left|\bigcup_{t=1}^T
(\tilde A_t\cap Q_k)\right|\\
= 1-\big(1-\mes_*(\tilde A_t) \big)^T> 1-\bigg(1-\frac{\mes_*(A_t)
}{32}\bigg)^T.
\end{multline*}
\end{proof}

For any line $l\subset \ZR^2$ we denote by $\arg l$ the positive
value of the minimal angle between $l$ and $x$-axes. For two
points $\theta ,\,\theta'\in\ZR^2$ we denote by $\theta \theta'$
the line passing through $\theta $ and $\theta'$, and by $[\theta,
\theta']$ the line segment with vertices $\theta$ and $\theta'$.
\begin{lemma}\label{lmR}
Let $0<\varepsilon <1$, $0<\gamma< \frac{\pi}{12}$ be any numbers
and
\begin{equation}\label{theta}
\theta_k=(\varepsilon/2^k,\,\sign (k)\cdot\tg\gamma\cdot
\varepsilon/2^k),\, k=\pm 1,\pm 2,\cdots .
\end{equation}
Then for any rectangle $R\in\MR_s$, with $3\gamma <|s|<
\frac{\pi}{2}-3\gamma $, we have
\begin{equation}\label{sumR}
\left|\sum_{0< |k|\le m,\,\theta_k\in R}\sign(k)\right|\le 2,\quad
m=1,2,\cdots .
\end{equation}
\end{lemma}
\begin{proof}
First we note that if $l$ is a line on the plane, then
\begin{equation}\label{geo}
l\cap [\theta_k,\theta_{-k}]\neq\varnothing ,\,l\cap
[\theta_{k+1},\theta_{-(k+1)}]\neq\varnothing
\end{equation}
implies
\begin{equation*}
 \arg l<3\gamma .
\end{equation*}

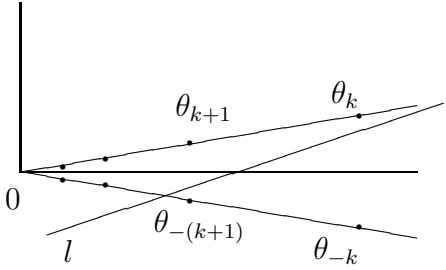
\begin{picture}(128,128)
\put(0,64){\line(1,0){150}} \put(0,64){\line(0,1){64}}
\put(0,64){\line(6,1){150}} \put(0,64){\line(6,-1){150}}
\put(10,40){\line(3,1){150}} \put(128,85){\circle*{2}}
\put(128,43){\circle*{2}} \put(64,75){\circle*{2}}
\put(64,53){\circle*{2}} \put(32,69){\circle*{2}}
\put(32,59){\circle*{2}} \put(16,66){\circle*{2}}
\put(16,61){\circle*{2}} \put(128,90){\llap{\sf $\theta_k$}}
\put(128,31){\llap{\sf $\theta_{-k}$}} \put(80,85){\llap{\sf
$\theta_{k+1}$}} \put(85,41){\llap{\sf $\theta_{-(k+1)}$}}
\put(0,50){\llap{\sf $0$}} \put(20,30){\llap{\sf $l$}}
\end{picture}

\noindent Indeed, using a simple geometry, one can check that
$\arg(\theta_{-k}\theta_{k+1})< 3\gamma$. Hence we get $\arg
l\le\arg(\theta_{-k}\theta_{k+1})<3\gamma $. Now consider a
rectangle
\begin{equation}\label{Rgam}
R\in\MR_s,\quad 3\gamma <|s|< \frac{\pi}{2}- 3\gamma.
\end{equation}
Let us show that
\begin{equation}\label{24}
\hbox { if }\theta_n,\theta_{n+1},\theta_{n+2}\in R,\hbox { then
}\theta_{-(n+1)}\in R.
\end{equation}
 Suppose we have the converse $\theta_{-(n+1)}\not\in R $. Then we can determine a
 line $l$ containing a side of $R$ and
separating the points $\theta_n,\theta_{n+1},\theta_{n+2}$ from
$\theta_{-(n+1)}$. Obviously we shall have
\begin{equation*}
l\cap [\theta_{n+1},\theta_{-(n+1)}]\neq\varnothing ,
\end{equation*}
and one of two following relations: $l\cap
[\theta_n,\theta_{-n}]\neq\varnothing $ or $l\cap
[\theta_{n+2},\theta_{-(n+2)}]\neq\varnothing $. So we have
\e{geo} for $k=n$ or $n+1$ and therefore $\arg l<3\gamma $, which
is a contradiction with \e{Rgam}. Similarly
\begin{equation}\label{25}
\hbox { if }\theta_{-n},\theta_{-(n+1)},\theta_{-(n+2)}\in R,\hbox
{ then }\theta_{n+1}\in R.
\end{equation}
Now let $p$ and $q$ are the numbers of elements of the sets
$\{1\le k\le m: \theta_k\in R\}$ and $\{-m\le k\le -1: \theta_k\in
R\}$. From \e{24} and \e{25} we conclude $|p-q|\le 2$, which
implies \e{sumR}.

\end{proof}

\begin{lemma}\label{lm2} For any numbers $0<\varepsilon <1$ and
$0<\gamma\le \frac{\pi}{12}$ there exists a bounded function $\phi
(x)=\phi (x_1,x_2)$ defined on $\ZR^2$ such that
\begin{gather}
\supp \phi\subset B(\varepsilon )
,\quad\int_{\ZR^2}\phi(x)dx=0,\quad\int_{\ZR^2}|\phi(x)|dx\le
1,\label{L31}
\\
\int_{\rot_s\big([0,x_1]\times [0,x_2
]\big)}\phi(x)dx\ge\frac{1}{4},\hbox { if  }\quad  x_1,x_2\ge
\varepsilon ,\quad |s|\le \gamma,\label{L32}
\\
M_s\phi (x)<\varepsilon , \hbox { as } x\not\in
\Gamma_s(2\varepsilon )\cup \Gamma_{s^\bot}(2\varepsilon )
,\,3\gamma <|s|< \frac{\pi}{2}- 3\gamma.\label{L33}
\end{gather}
\end{lemma}
\begin{proof}
Consider the sequence $\theta=\theta^+\cup\theta^-$ where
\begin{multline}\label{thetaN}
\theta^+=\{\theta_k:\quad k=1,2,\cdots ,N\},\,\\
\theta^-=\{\theta_k:\quad k=-1,-2,\cdots ,-N\},\
N=[10\varepsilon^{-3} ]+1,
\end{multline}
and $\theta_k$ are defined in \e{theta}. We have
\begin{equation*}
\theta_k\in B\bigg(\frac{\varepsilon}{\sqrt 2}\bigg)\subset
B(\varepsilon ), \quad \theta_k\in\{x:x_2=\tg \gamma \cdot x_1 \}
\quad k=\pm 1,\pm 2,\cdots .
\end{equation*}
Define the balls $b_k$, denoting
\begin{equation*}
 b_k=\{x\in R^2: |x-\theta_k|< r\},\quad
 k=\pm 1,\pm 2,\cdots ,\pm N.
\end{equation*}
Choosing a small number $r>0$, we provide the following
conditions:
\begin{itemize}
    \item [1)] $b_k\subset B(\varepsilon )$ and they are mutually disjoint,
    \item [2)] if $k>0$, then $b_k$ is in the upper half-plane, if
    $k<0$ is in lower,
    \item [3)]any line $l$ with $|\arg l|\ge 3\gamma $ intersects
at most two $b_k$.

\end{itemize}
 We define
\begin{equation*}
\phi(x)=\frac{1}{2\pi Nr^2}\sum_{k=1}^N
\big(\ZI_{b_k}(x)+\ZI_{b_{-k}}(x)\big),
\end{equation*}
where $\ZI_{b_k}$ is the characteristic function of $b_k$. The
conditions \e{L31} are clear. To show \e{L32} we shall use
conditions 1) and 2). We fix numbers $x_1,x_2>\varepsilon $. If
$0\le s<\gamma $, then we have
\begin{gather*}
\rot_s\big([0,x_1]\times [0,x_1]\big)\cap b_k=\varnothing \hbox {
as } -N\le k<0,
\\
|\rot_s\big([0,x_1 ]\times [0,x_2]\big)\cap
b_k|>\frac{|b_k|}{2}=\frac{\pi r^2}{2}\hbox { as }0< k\le N.
\end{gather*}
Therefore
\begin{equation*}
\int_{\rot_s\big([0,x_1 ]\times [0,x_2
]\big)}\phi(x)dx=\frac{1}{2\pi
Nr^2}\sum_{k=1}^N\int_{\rot_s\big([0,x_1 ]\times [0,x_2
]\big)}\ZI_{b_k}(x)dx\ge  \frac{1}{4}.
\end{equation*}
If $-\gamma < s\le 0 $, then
\begin{gather*}
b_k\subset\rot_s\big([0,x_1 ]\times [0,x_1 ]\big),\quad k>0,
\\
|\rot_s\big([0,x_1 ]\times [0,x_2 ]\big)\cap b_k|\le
\frac{|b_k|}{2}=\frac{\pi r^2}{2},\quad k>0,
\end{gather*}
and then similarly we obtain \e{L32}. We shall prove now if
\begin{equation}\label{26}
R\in\MR_s,\,3\gamma <|s|< \frac{\pi}{2}- 3\gamma
\end{equation}
then
\begin{equation}\label{intphi}
\left|\int_R\phi (x)dx\right|\le \frac{10}{N}<\varepsilon^3.
\end{equation}
We have
\begin{multline}\label{sum}
\int_R\phi (x)dx=\frac{1}{2\pi Nr^2} \sum_{b_k\cap
R\neq\varnothing } \int_R\ZI_{b_k}(x)dx=
\\
\frac{1}{2\pi Nr^2}\sum_{\theta_k\in R}
\int_R\ZI_{b_k}(x)dx+\frac{1}{2\pi Nr^2}\sum_{\theta_k\not\in
R,b_k\cap R\neq\varnothing } \int_R\ZI_{b_k}(x)dx.
\end{multline}
The conditions $\theta_k\not\in R,b_k\cap R\neq\varnothing $ mean
that $b_k$ intersects a side of $R$. Also we have that if a line
$l$ contains a side of $R$ then $|\arg l|>3\gamma $. On the other
hand by the condition 3) any line with $|\arg l|>3\gamma $ can
intersect not more than two balls $b_k$. So the number of terms in
the second sum doesn't exceed $8$. Therefore
\begin{equation}\label{sum2}
\left|\frac{1}{2\pi Nr^2}\sum_{\theta_k\not\in R,b_k\cap
R\neq\varnothing } \int_R\ZI_{b_k}(x)dx\right|\le \frac{4}{N}.
\end{equation}
By the same reason the equality
\begin{equation*}
\int_R\ZI_{b_k}(x)dx=\int_{\ZR^2}\ZI_{b_k}(x)dx
\end{equation*}
fails for not more than $8$ different $k$'s. Therefore
\begin{equation*}
\left|\frac{1}{2\pi Nr^2}\sum_{\theta_k\in R}
\int_R\ZI_{b_k}(x)dx-\frac{1}{2\pi Nr^2}\sum_{\theta_k\in R}
\int_{\ZR^2}\ZI_{b_k}(x)dx\right|\le\frac{4}{N}.
\end{equation*}
Hence we obtain
\begin{multline}\label{sum1}
\left|\frac{1}{2\pi Nr^2}\sum_{\theta_k\in R}
\int_R\ZI_{b_k}(x)dx\right|\le \left|\frac{1}{2\pi
Nr^2}\sum_{\theta_k\in R}
\int_{\ZR^2}\ZI_{b_k}(x)dx\right|+\frac{4}{N}=
\\
\left|\frac{1}{2 N} \sum_{\theta_k\in R}
\sign(k)\right|+\frac{4}{N}\le \frac{5}{N},
\end{multline}
where the last inequality follows from the \lem{lmR}. Combining
\e{sum}, \e{sum1} and \e{sum2} we get \e{intphi}. Fix a slope $s$
with $3\gamma <|s|\le \frac{\pi}{4}$ and take a point $x\in \ZR^2$
such that
\begin{gather*}
x\not\in \Gamma_s(2\varepsilon )\cup \Gamma_{s^\bot}(2\varepsilon
),
\\
 x\in R\in\MR_s,\,3\gamma <|s|< \frac{\pi}{2}- 3\gamma.
\end{gather*}
We need to prove
\begin{equation}\label{Rvar}
\frac{1}{|R|}\int_R\phi(t)dt\le \varepsilon .
\end{equation}
Assume the lengths of the sides of $R$ are $a$ and $b$. If $R$
doesn't contain a point $\theta_k$ then \e{Rvar} is trivial. So we
suppose there exists at least one point $\theta_k\in R$. Hence $R$
has an intersection with  $B(\varepsilon )$ and
$\big(\Gamma_s(2\varepsilon )\cup \Gamma_{s^\bot }(2\varepsilon
)\big)^c$. Taking account of $R\in \MR_s$ we get
$a,b>\varepsilon$. Hence by \e{intphi} we get
\begin{equation*}
\frac{1}{|R|}\int_R\phi(t)dt\le \frac{\varepsilon^3}{ab}\le
\varepsilon
\end{equation*}
\end{proof}

\begin{lemma}\label{lm3}For any numbers  $0<\varepsilon,\delta  <1/10$,
and interval $S=[\alpha-\gamma,\alpha+\gamma ]\subset [0,\pi/2 )$
with $0<\gamma \le\frac{\pi }{12}$ there exist a bounded function
$\phi (x)$ and numbers $\nu ,\nu'$ with $0<\nu <\nu'$ such that
\begin{gather}
\sup_{k\in \ZZ^2} \int_{Q_k}|\phi(x)|dx\le 1\label{l1}
\\
\mes^*\{x\in \ZR^2:\, M_s\phi (x)>\varepsilon \}<\varepsilon,\quad
3\gamma <|s-\alpha |< \frac{\pi}{2}- 3\gamma ,\label{l3}
\\
\mes^*\{x\in \ZR^2:\, M_s^{[0,\nu )}\phi (x)>\varepsilon
\}<\varepsilon,\quad s\in [0,2\pi ),\label{l5}
\\
 M_s^{[\nu' ,\infty)}\phi (x)<\varepsilon ,\quad x\in \ZR^2,\,
s\in [0,2\pi ),\label{l4}
\\
 \mes_*\{
M_s^{[\nu, \nu'] }\phi (x)>
\frac{1}{\delta}\}>\frac{\delta}{4}\ln\frac{1}{12\delta},\quad
s\in S.\label{l2}
\end{gather}
\end{lemma}
\begin{proof} Without loss of generality we may
assume $\alpha=0$, i.e. $S=[-\gamma,\gamma ]$. We take $\lambda
=\min \{\varepsilon /100,\delta \}$ and consider a double sequence
$\varepsilon_k=\varepsilon_{k_1,k_2} =\lambda
2^{-(|k_1|+|k_2|)},\, k\in \ZZ^2$. Using \lem{lm2} we can find
functions $\phi_k(x) $ with following conditions:
\begin{gather}
\supp \phi_k\subset B(\varepsilon_k ) \subset B(\varepsilon ),
\label{sup}
\\
\int_{Q_0}\phi_k(x)dx=0,\quad\int_{Q_0}|\phi_k(x)|dx\le 1,
\\
\int_{\rot_s\big(R_x\big)}\phi_k(x)dx>\frac{1}{4},\, R_x=[0,x_1
]\times [0,x_2 ],\, x_1,x_2\ge\delta \ge\varepsilon_k ,\,
|s|<\gamma,\label{integ}
\\
M_s\phi_k (x)<\varepsilon_k , \hbox { as } x\not\in
\Gamma_s(2\varepsilon_k )\cup \Gamma_{s^\bot}(2\varepsilon_k )
,\,3\gamma <|s|\le \frac{\pi}{2}-3\gamma,\label{Msphi}
\end{gather}
where $k=(k_1,k_2)$. Denote
\begin{gather}
\phi(x)=\sum_{k\in \ZZ^2}\phi_k(x+k),\label{defphi}
\\
E_s=\bigcup_{k\in \ZZ^2}\bigg(k+\big(\Gamma_s(2\varepsilon_k )\cup
\Gamma_{s^\bot}(2\varepsilon_k )\big)\bigg).\label{defE}
\end{gather}
We obviously have \e{l1} and
\begin{gather}
\supp \phi (x)\subset \bigcup_{k\in \ZZ^2}\big(k+B(\varepsilon
)\big),\label{phi1}
\\
\int_{Q_k}\phi(x)dx=0,\quad k\in \ZZ^2. \label{phi2}
\end{gather}

\underline {Proof of \e{l3}}: For any square $Q_j$, $j\in \ZZ^2$,
we have
\begin{gather*}
|Q_j\cap \big(k+\Gamma_s(2\varepsilon_k )\big)|\le \diam
Q_j\times\mes_s\big(k+\Gamma_s(2\varepsilon_k )\big)=
4\varepsilon_k\sqrt{2},
\\
|Q_j\cap\big(k+\Gamma_{s^\bot}(2\varepsilon_k )\big)|\le
4\varepsilon_k\sqrt{2}.
\end{gather*}
Hence we obtain
\begin{equation}\label{Es}
\mes^*(E_s)\le\sum_k 8\sqrt{2}\varepsilon_k=32\sqrt{2}\lambda \le
\varepsilon .
\end{equation}
From \e{Msphi} it follows that
\begin{equation*}
M_s\phi_k(x+k)\le \varepsilon_k,\quad x\not\in E_s\supset
k+\big(\Gamma_s(2\varepsilon_k )\cup
\Gamma_{s^\bot}(2\varepsilon_k )\big),\quad 3\gamma <|s|\le
\frac{\pi}{2}-3\gamma.
\end{equation*}
 Then according \e{defphi} and \e{defE} we get
\begin{equation*}
M_s\phi(x)\le \sum_kM_s\phi_k(x+k)\le
\sum_k\varepsilon_k\le\varepsilon, x\not\in E_s,\quad 3\gamma
<|s|\le \frac{\pi}{2}-3\gamma,
\end{equation*}
and combining this with \e{Es} we obtain \e{l3}.

\underline{Proof of \e{l5}}: From \e{phi1} it follows that
\begin{equation*}
\lim_{\nu\to 0}M_s^{[0,\nu )}\phi (x)=0,\, \hbox{ if }x\not\in
\bigcup_{k\in \ZZ^2}\big(k+B(\varepsilon )\big),\quad s\in
[0,2\pi),
\end{equation*}
therefore for a small $\nu <\delta $ we shall have \e{l5}, since
\begin{equation*}
\mes^*\bigg(\bigcup_{k\in \ZZ^2}\big(k+B(\varepsilon
)\big)\bigg)=|B(\varepsilon )|=\pi \varepsilon^2\le \varepsilon .
\end{equation*}
\underline{Proof of \e{l4}}: From \e{phi2} we obtain
\begin{equation*}
\lim_{\nu'\to\infty }\int_R\phi(\nu' x)dx=0
\end{equation*}
for any rectangle $R$ and the convergence is uniformly by $R\in
\MR_s^{[1,\infty )},\, s\in [0,2\pi )$. So for a big $\nu'>1/4 $
we shall have
\begin{equation*}
M_s^{[1 ,\infty)}\phi (\nu' x)<\varepsilon ,\,x\in \ZR^2,\quad
s\in [0,2\pi).
\end{equation*}
By dilation we get
\begin{equation*}
M_s^{[\nu' ,\infty)}\phi (x)=M_s^{[1 ,\infty)}\phi (\nu'
x)<\varepsilon ,\,x\in \ZR^2,\quad s\in [0,2\pi),
\end{equation*}
which gives \e{l4}.

\underline{Proof of \e{l2}}: Consider the set
\begin{equation}\label{Adelta}
A=\big\{x=(x_1,x_2):\, x_1x_2\le\frac{\delta}{4} ,\quad \delta \le
x_1,x_2\le\frac{1}{4}\big\}.
\end{equation}
We have
\begin{gather}
\rot_sA\subset B\bigg(\frac{1}{2}\bigg),\quad s\in
[-\pi/4,\pi/4),\label{Ames1}\\
|A|=\int_{\delta}^{1/4} \frac{\delta}{4t}dt-\delta
\bigg(\frac{1}{4}-\delta
\bigg)>\frac{\delta}{4}\ln\frac{1}{12\delta}\label{Ames2}
\end{gather}
If $x=(x_1,x_2)\in A$, then
\begin{equation}\label{x1x2}
\frac{1}{4}\ge x_1,x_2\ge \delta>\varepsilon_k ,\quad |R_x|\le
\frac{\delta}{4}
\end{equation}
 So by \e{integ}
\begin{equation*}
\int_{k+\rot_s(R_x)}\phi_k(k+t)dt =\int_{\rot_s(R_x)}\phi_k(t)dt
>\frac{1}{4},\hbox { as } x\in A,\,|s|<\gamma ,
\end{equation*}
and therefore from \e{sup} we can get
\begin{equation}\label{krot}
\int_{k+\rot_s(R_x)}\phi(t)dt=\int_{k+\rot_s(R_x)}\phi_k(k+t)dt>
\frac{1}{4},\hbox { as }\,x\in A,\, |s|< \gamma.
\end{equation}
According to $\nu <\delta $, $\nu'>1/4$ we have $R_x\in
\MR_0^{[\delta,1/4 ]}\subset\MR_0^{[\nu,\nu' ]}$. Since $|R_x|\le
\delta/4 $  from \e{krot} and \e{x1x2} we conclude
\begin{equation}\label{MGs}
M_s^{[\nu,\nu' ]}\phi(x)>\frac{1}{4|R_x|}>\frac{1}{\delta },\quad
x\in G_s=\bigcup_k\big(k+\rot_s A \big),\quad |s|< \gamma.
\end{equation}
In addition, by \e{Adelta}, \e{Ames1} and \e{Ames2}, for any $m\in
\ZZ^2$ we get
\begin{equation*}
|(m+Q_0)\cap G_s |=|m+\rot_s A|=|A|
>\frac{\delta}{4}\ln\frac{1}{12\delta},
\end{equation*}
which implies
\begin{equation*}
\mes^*G_s>\frac{\delta}{4}\ln\frac{1}{12\delta}.
\end{equation*}
Combining this with \e{MGs} we obtain \e{l2}.
\end{proof}
\end{section}

\begin{section}{Proofs of Theorems}

\begin{proof}[Proof of Theorem 1]
Let $G$ be an arbitrary $G_\delta$-set in $[0,\pi/2)$. So
\begin{equation*}
G=\big(\cap_{k=1}^\infty G_k\big)\cap [0,\pi/2),
\end{equation*}
where $G_k\subset \ZR $ are open sets and
\begin{equation*}
G_1\supseteq G_2\supseteq \cdots \supseteq G_n\supseteq\cdots .
\end{equation*}
Each $G_k$ is union of a mutually disjoint intervals, i.e.
\begin{equation*}
G_k=\cup_jI_j^k.
\end{equation*}
We note that an arbitrary interval $I=(\alpha ,\beta )\subset \ZR
$ can be split to disjoint intervals $I_i=[\alpha_i,\beta_i)$ such
that
\begin{equation*}
|I_i|\le \frac{\pi}{12},\quad 3I_i\subset I,\quad
\sum_i\ZI_{3I_i}(x)\le 8.
\end{equation*}
For $I=(-1,1)$ such a partition is
\begin{gather*}
\bigg[1-\bigg(\frac{9}{10}\bigg)^k,1-
\bigg(\frac{9}{10}\bigg)^{k+1}\bigg),\quad k=0,1,2,\cdots,
\\
\bigg[\bigg(\frac{9}{10}\bigg)^{k+1}-1,
\bigg(\frac{9}{10}\bigg)^k-1\bigg),\quad k=0,1,2,\cdots,
\end{gather*}
We do a similar splitting for any $I_j^k$. Let $J_t,\,
t=1,2,\cdots ,$ be a numeration of those splitting intervals $J$
for wich $J\cap [0,\pi/2)\neq \varnothing $. We denote
$l_t=J_t\cap [0,\pi/2)$. It is easy to check the following two
relations
\begin{itemize}
    \item [1)] if $x\in G$, then $x$ belongs to infinite number of
    $l_t$'s,
    \item [2)] if $x\not\in G$ then $x$ belongs only to finite number of
    $3l_t$'s.
\end{itemize}
We chose integers $0=m_0<m_1<m_2<\cdots $ satisfying
\begin{equation}\label{33}
\prod_{k=m_t+1}^{m_{t+1}}\bigg(1-\frac{1}{k\ln
k}\bigg)<\frac{1}{2^t},\quad t=1,2,\cdots .
\end{equation}
We denote
\begin{equation}
S_k=l_t,\hbox{ if } m_t<k\le m_{t+1}.
\end{equation}
 Using \lem{lm3} for $S=S_k$, $\varepsilon =1/2^k$,
 $\delta =1/k\ln^2k$, we may
define functions $\phi_k(x)$ and numbers $0<\nu_k<\nu'_k$ with
conditions \e{l1}-\e{l2}. We denote
\begin{gather}
U_{s,k}=\{x\in \ZR^2:\, M_s\phi_k (x)\le\frac{1}{2^k}\},\label{c1}\\
V'_{s,k}=\{x\in \ZR^2:\, M_s^{[0,\nu_k)}\phi_k (x)
\le\frac{1}{2^k}
\},\label{c2}\\
V''_{s,k}=\{ x\in \ZR^2:M_s^{[\nu_k , \nu'_k ] }\phi_k (x)>
k\ln^2k\}.\label{c3}
\end{gather}
By \e{l3},\e{l5},\e{l2} we have
\begin{gather}
\mes_*U_{s,k}>1-\frac{1}{2^k},\quad s\in [0,\pi/2 )\setminus
3S_k,\label{e1}
\end{gather}
(we may replace the condition $3\gamma <|s-\alpha |<
\frac{\pi}{2}- 3\gamma$ in \e{l3} by  $s\in [0,\pi/2 )\setminus
3S$ because the second implies the first) and
\begin{gather}
\mes_*V'_{s,k}>1-\frac{1}{2^k},\quad s\in [0,\pi/2 ),\label{e2}\\
 \mes_*V''_{s,k}>\frac{1}{4k\ln^2k}\ln\frac{k\ln^2k}{12}>\frac{c}{k\ln
k},\quad s\in S_k\quad (k\ge 3).\label{e3}
\end{gather}
From \e{l4} we get
\begin{equation}
M_s^{[\nu'_k ,\infty)}\phi_k (x)<\frac{1}{2^k},\quad x\in \ZR^2,\,
s\in [0,\pi/2 ).
\end{equation}
We define integers $1=n_0<n_1<n_2<\cdots $, so that
\begin{equation}\label{34}
\frac{n_k}{n_{k-1}}>\max\bigg(\frac{4}{\nu_{k-1}},
\frac{\nu'_k}{\nu_{k-1}}\bigg),\quad k=1,2,\quad ,
\end{equation}
and denote $\mu_k=\nu_k/n_k$. It is clear
\begin{equation*}
\mu_{k-1}>\frac{\nu'_k}{n_k}>\mu_k,\quad k=2,3,\cdots .
\end{equation*}
Consider the functions
\begin{equation}\label{310}
\psi_k(x)=\phi_k(n_kx),\quad x\in Q_0.
\end{equation}
According to \e{c1}-\e{c3} and \e{310}, we obviously have
\begin{gather}
M_s\psi_k (x)\le\frac{1}{2^k},\quad x\in\dil_{n_k}U_{s,k},
\quad s\in [0,\pi/2 )\setminus 3S_k,\label{b1}\\
 M_s^{[0,\mu_k)}\psi_k (x)= M_s^{[0,\nu_k/n_k)}\psi_k (x) \le\frac{1}{2^k}\quad x\in\dil_{n_k}V'_{s,k},
 \quad s\in [0,\pi/2 ),\label{b2}\\
M_s^{[\mu_k , \mu_{k-1}] }\psi_k (x)> M_s^{[\nu_k/n_k , \nu'_k/n_k
] }\psi_k (x)> k\ln^2k,\quad
 x\in\dil_{n_k}V''_{s,k},\quad s\in S_k,\label{b3}\\
 M_s^{[\mu_{k-1},\infty)}\psi_k (x)\le M_s^{[\nu'_k/n_k,\infty)}\psi_k (x)
  \le\frac{1}{2^k},\quad x\in\ZR^2,
 \quad s\in [0,\pi/2 )\label{b4}.
\end{gather}
Desired function will be
\begin{equation}\label{32}
f(x)=\sum_{k=1}^\infty \frac{\psi_k (x)}{k\ln^{3/2}k},\quad x\in
Q_0.
\end{equation}
Denote
\begin{equation}\label{36}
U_s= \limsup_{k\to\infty }\bigg(\big(\dil_{n_k}U_{s,k}\big)\cap
Q_0\bigg),
\end{equation}
where $\limsup_{k\to\infty }A_k$ means $\cup_n\cap_{k\ge n}A_k$.
 If $s\not\in G$, then by 2) $s\in [0,\pi/2 )\setminus 3S_k$ as $k>k(s)$.
Therefore, by \e{e1} we have $|\dil_{n_k}U_{s,k}\cap
Q_0|\ge\mes_*U_{s,k}
>1-1/2^k,\, k>k(s),$ and so we get
\begin{equation}
|U_s|=1\hbox { if } s\not\in G.
\end{equation}
From \e{b1} and \e{36} we get, that for any $x\in U_s$
\begin{equation*}
M_s\psi_k(x)\le\frac{1}{2^k},\quad k>k(x).
\end{equation*}
Hence, if $\varepsilon >0$, then for an appropriate $N>k(x)$ we
get
\begin{equation}\label{31}
M_s\bigg(\sum_{k=N+1}^\infty
\frac{\psi_k(x)}{k\ln^{3/2}k}\bigg)\le \sum_{k=N+1}^\infty
\frac{M_s\psi_k(x)}{k\ln^{3/2}k}\le
\sum_{k=N+1}^\infty\frac{1}{k2^k\ln^{3/2}k}<\varepsilon  .
\end{equation}
On the other hand, since
\begin{equation*}
\sum_{k=1}^N \frac{\psi_k(x)}{k\ln^{3/2}k}
\end{equation*}
is a bounded function, the basis $\MR_s$ differentiates its
integral. So, taking account of \e{31} and \e{32} we get $\int f$
differentiable by $\MR_s$ if $s\in [0,\pi/2)\setminus G$.

Now let us take $s\in G$. We have $s\in l_{t_i}$, $i=1,2,\cdots $.
Hence $s\in S_k$ if $m_{t_i}<k\le m_{t_i+1},\, i=1,2,\cdots $. We
notice, that each $V''_{s,k}$ defined in \e{c3} is $\nu_k$-set,
and by \e{34} $n_{k+1}>\frac{4}{\nu_k}n_k$. Therefore, using
\e{33},
 from \lem{lm5} we obtain
\begin{equation}\label{37}
\left|\bigcup_{k=m_{t_i}+1}^{m_{t_i+1}}\dil_{n_k}V''_{s,k}\cap
Q_0\right|\ge
1-\prod_{k=m_{t_i}+1}^{m_{t_i+1}}\bigg(1-\frac{1}{k\ln
k}\bigg)>1-\frac{1}{2^t}.
\end{equation}
Denoting
\begin{equation*}
V_s=\bigg(\limsup_{k\to\infty
}\dil_{n_k}V'_{s,k}\bigg)\bigcap\bigg(\limsup_{i\to\infty }
\bigcup_{k=m_{t_i}+1}^{m_{t_i+1}}\dil_{n_k}V''_{s,k}\bigg)\bigcap
Q_0,
\end{equation*}
from \e{37} and \e{e2} we get
\begin{equation}\label{38}
|V_s|=1,\quad s\in G.
\end{equation}
On the other hand if $x\in V_s$, then
\begin{gather*}
x\in \dil_{n_{k_i}}V''_{s,k_i},\quad i=1,2,\cdots ,\\
x\in\dil_{n_k}V'_{s,k},\quad k>k(x).
\end{gather*}
 where $k_i\to \infty $, and therefore, by \e{b2} and
 \e{b4}  we have
\begin{equation*}
M_s^{[\mu_{k_i} , \mu_{{k_i}-1}] }\psi_j (x)\le
\frac{1}{2^{k_i}},\hbox { if } j\neq k_i.
\end{equation*}
The case $j>k_i$ follows from \e{b4} and $j<k_i$ from \e{b2}. From
\e{b3} we get
\begin{equation*}
 M_s^{[\mu_{k_i} , \mu_{k_i-1}] }\psi_{k_i} (x)> k_i\ln^2k_i
.
\end{equation*}
 So if $k_i>k(x)$, then
\begin{multline*}
M_sf(x)\ge M_s^{[\mu_{k_i} , \mu_{k_i-1}]}f(x)\ge
\\
\frac{M_s^{[\mu_{k_i},
\mu_{k_i-1}]}\psi_{k_i}(x)}{k_i\ln^{3/2}k_i}-\sum_{j\neq k_i}
\frac{M_s^{[\mu_{k_i} , \mu_{k_i-1}]}\psi_j(x)}{j\ln^{3/2}j}\ge
c\sqrt{\ln k_i}-\sum_{j\neq k_i} \frac{1}{j2^j\ln^{3/2}j}
\end{multline*}
and so $\overline D_sf(x)=\infty $, whenever $x\in V_s$ and $s\in
G$. Since $|V_s|=1$ by \e{38}, the theorem is completely proved.
\end{proof}
\begin{proof}[Proof of Theorem 2] The necessity of
the theorem is shown in the introduction. To prove the sufficiency
we let $V\in[0,\pi/2)$ to be an arbitrary $G_{\delta\sigma }$ set
and we have
\begin{equation*}
V=\bigcup_n V_n
\end{equation*}
where each $V_n$ is $G_\delta $. According to Theorem 1 for each
$V_n$ there exists a function $f_n\in L^1(\ZR^2)$ such that its
integral differentiable by $\MR_s$ as $s\not\in  V_n$ and
$\overline D_s f_n(x)=\infty $ a.e. if $s\in  V_n$. Denote
$g_n(x)=\chi_{Q_n}(x)f_n(x)$, where $Q_n$ is a family of arbitrary
pairwise disjoint unit open squares, and consider the function
\begin{equation*}
g(x)=\sum_{n=1}^\infty g_n(x).
\end{equation*}
 Since the supports
of the functions $g_n$ are disjoint for any point $x\in Q_n$ and
any $s$ we have
\begin{equation*}
\overline D_s g(x)= \overline D_s g_{n}(x)=\overline D_s f_{n}(x).
\end{equation*}
If $s\in V$ then $s\in V_n$ for some $n$. So we get $\overline
D_sg(x) =\overline D_s f_{n}(x)=\infty $ almost everywhere on the
square $Q_n$. Using disjointness of the supports of the functions
$g_n$ once again we conclude that if $s\not\in V$ then
\begin{equation*}
\lim\limits_{d(R)\to 0, x\in R\in\MR_s} \frac{1}{|R|}\int_R
g=g(x)\hbox { a.e. }.
\end{equation*}
Finally we get that $V$ is $WR$-set and Theorem 2 is proved.
\end{proof}

\end{section}

\end{document}